\documentclass[reqno]{amsart}
\RequirePackage[l2tabu, orthodox]{nag}
\usepackage{amssymb} 
\usepackage{amscd}
\usepackage{bbm}
\usepackage{color}
\usepackage{esint}
\usepackage{framed}
\usepackage{graphicx}
\usepackage{hyphenat}
\usepackage{mathrsfs} 
\usepackage{lmodern}
\usepackage{microtype}
\usepackage{hyperref} 
\usepackage{amsrefs}

\makeatletter
\DeclareFontFamily{OMX}{MnSymbolE}{}
\DeclareSymbolFont{MnLargeSymbols}{OMX}{MnSymbolE}{m}{n}
\SetSymbolFont{MnLargeSymbols}{bold}{OMX}{MnSymbolE}{b}{n}
\DeclareFontShape{OMX}{MnSymbolE}{m}{n}{
    <-6>  MnSymbolE5
   <6-7>  MnSymbolE6
   <7-8>  MnSymbolE7
   <8-9>  MnSymbolE8
   <9-10> MnSymbolE9
  <10-12> MnSymbolE10
  <12->   MnSymbolE12
}{}
\DeclareFontShape{OMX}{MnSymbolE}{b}{n}{
    <-6>  MnSymbolE-Bold5
   <6-7>  MnSymbolE-Bold6
   <7-8>  MnSymbolE-Bold7
   <8-9>  MnSymbolE-Bold8
   <9-10> MnSymbolE-Bold9
  <10-12> MnSymbolE-Bold10
  <12->   MnSymbolE-Bold12
}{}

\let\llangle\@undefined
\let\rrangle\@undefined
\DeclareMathDelimiter{\llangle}{\mathopen}%
                     {MnLargeSymbols}{'164}{MnLargeSymbols}{'164}
\DeclareMathDelimiter{\rrangle}{\mathclose}%
                     {MnLargeSymbols}{'171}{MnLargeSymbols}{'171}
\makeatother

\theoremstyle{plain} 
\newtheorem{theorem}{Theorem}[section]
\newtheorem{proposition}[theorem]{Proposition}
\newtheorem{lemma}[theorem]{Lemma} 

\theoremstyle{definition} 
\newtheorem{definition}[theorem]{Definition}

\theoremstyle{remark} 
\newtheorem{remark}[theorem]{Remark}

\definecolor{shadecolor}{rgb}{1,0.8,0.3}
\newcommand{\DETAIL}[1]{}
\newcommand{\IGNORE}[1]{}

\newcommand{\N}{{\mathbf{N}}}
\newcommand{\R}{{\mathbf{R}}}

\newcommand{\ID}{{\mathrm{id}}}
\newcommand{\LS}{\leqslant}

\newcommand{\PR}{{\SP\SUB{\RHO}(\R^{2d})}}

\newcommand{\SP}{{\mathscr{P}}}
\newcommand{\TF}{{\mathcal{T}}}

\newcommand{\WR}{{\WAS\SUB{\RHO}}}

\newcommand{\ADM}{{\mathrm{ADM}}}

\newcommand{\ETA}{{\boldsymbol{\eta}}}

\newcommand{\OPT}{{\mathrm{OPT}}}

\newcommand{\RHO}{\varrho}
\newcommand{\SUB}[1]{_{\raisebox{0.25ex}{\scriptsize{$#1$}}}}

\newcommand{\WAS}{{\mathrm{W}}}
\newcommand{\BETA}{{\boldsymbol{\beta}}}
\newcommand{\COPT}{C}

\newcommand{\ALPHA}{{\boldsymbol{\alpha}}}
\newcommand{\GAMMA}{{\boldsymbol{\gamma}}}

\DeclareMathOperator{\SPT}{{\mathrm{spt}}}

\numberwithin{equation}{section}

\title[Tangent Cone of Monotone Transport Plans]{Characterization of  the Tangent space of monotone transport plans in ${\R}\times{\R}$ with prescribed first projection.}

\author
  {Marc Sedjro}
\address
  {Marc Sedjro,
   Lehrstuhl f\"{u}r Mathematik (Analysis),
   RWTH Aachen University,
   Templergraben 55,
   D-52062 Aachen, 
   Germany}
\email
  {sedjro@instmath.rwth-aachen.de}

\author
  {Michael Westdickenberg}
\address
  {Michael Westdickenberg,
   Lehrstuhl f\"{u}r Mathematik (Analysis),
   RWTH Aachen University,
   Templergraben 55,
   D-52062 Aachen, 
   Germany}
\email
  {mwest@instmath.rwth-aachen.de}

\date{\today}
\subjclass[2000]{ 46B20, 51F99,}
\keywords{, tangent cone, monotone transport plans, Wasserstein metric optimal transport.}

\begin{document}

\begin{abstract}
We provide a complete description of  the tangent space of the cone of monotone plans in $\R\times \R$ with prescribed first projection. We show  that elements of this tangent space are essentially made of two simple building-block types of measures.  
\end{abstract}

\maketitle

 \section{Introduction.}\label{Def}
The study of  continuity equation has shown that an appropriate notion of  tangent space in the 2-Wasserstein space $\SP_{2}(\R^{d})$, at $\RHO\in \SP_{2}(\R^{d})$ is given by: 
\begin{equation}\label{tangent space cont eq} 
 \text{Tan}_{\RHO}\SP_{2}(\R^{d})=\overline{\left\lbrace \nabla\varphi :\varphi\in C_{c}(\mathbb{R^d})\right\rbrace }^{L^2(\RHO,\R^{d})} \qquad d\in \N.
\end{equation}
While this notion of tangent space has proven to be a valuable tool in many problems, it has also shown its limit especially when dealing with tangent space at non-regular measures.  For instance, when $\RHO$ is equal to the average of $n$ delta masses sitting at $n$ distinct points, $\text{Tan}_{\RHO}\SP_{2}(\R^{d})$ is isomorphic to $\R^{nd}$, which makes it  too small to be of great interest.

In \cite{AmbrosioGigliSavare2005} and \cite{Gigli2004}, the authors introduced a broader notion of tangent spaces called geometric tangent cone : 

\begin{equation}\label{ geom tangent cone } 
 \textbf{Tan}_{\RHO}\SP_{2}(\R^{d}):=\overline{\left\lbrace \gamma\in\SP_{\RHO}(\R^{2d}) :\left(\pi^1,\pi^1+\varepsilon\pi^2 \right)\#\gamma \text{  optimal for some   }\varepsilon >0\right\rbrace }^{  \WAS\SUB{\RHO}}
\end{equation}
Here, $\SP_{\RHO}(\R^{2d})$ denotes the set of transport plans whose first marginal is $\RHO$ and $\WAS\SUB{\RHO}$, its associated metric; see Definition \ref{def metric}. $\pi^1,\pi^2$ denote respectively, the first and second projection maps. It is important to note that $\textbf{Tan}_{\RHO}\SP_{2}(\R^{d})$ is constructed by using the geodesics with respect to the $2$-Wasserstein distance. Even though the two tangent spaces  come from  different perspectives, it was shown that $\text{Tan}_{\RHO}\SP_{2}(\R^{d})$ is isometrically embedded in  $\textbf{Tan}_{\RHO}\SP_{2}(\R^{d})$. As a consequence, the geometric tangent cone offers a richer supply of elements. Behind the well-behaved nature of   $\textbf{Tan}_{\RHO}\SP_{2}(\R^{d}) $ hides the fact that this cone adapts well to the possibility of mass splitting that is precisely missing in  (\ref{tangent space cont eq}).

In spite of the gain the tangent space concept in (\ref{ geom tangent cone })  offers over the former concept in (\ref{tangent space cont eq}), it is still far from being satisfactory.  In the recent work \cite{CavSedWes2014}, the  authors extended the Lagrangian formulation of the  one dimensional  pressureless Euler equations to the multidimensional case. It turns out that the configuration space is better described by \textit{monotone transport plans}, that is transport plans whose support is a monotone set. This means that relying on geodesics could be too restrictive when it comes to  describing the configuration space.\\

 In this work, we propose a slightly different concept of tangent cone from convex analysis point of view :  the tangent cone of the cone of  monotone transport plans  whose first projection is $\RHO$ in the space $\SP_{\RHO}(\R^{2d})$. It is noteworthy to mention that in dimension one,  the set of optimal maps coincides with the set of monotone maps. We initiate  the study of the tangent cone of  this cone of monotone transport plans and give a complete characterization  for $d=1$. In the framework developed in \cite{Gigli2004},  we show that this tangent cone is convex. As a result, elements of the tangent cone  are linear combinations of two basic types of measures so that on the diffuse part on $\RHO$, an element of the tangent cone is induced by a map and, at each atomic point its support can be any closed subsets of $\R$.
 
  This paper is organized in the following way:  section \ref{sec: notation} provides some notation, definitions, and general results related to $\SP_{\varrho}(\R^2)$  and its associated metric. In Section \ref{Prel}, we give some convergence results and state the main theorem of the paper. In Section \ref{Useful}, the tangent space of the monotone transport plans is proven to be a convex cone in $\SP_{\RHO}(\R^{2d})$. Section \ref{Chara of the Tangent} deals with  the complete description of the tangent space of monotone plans in $\SP_{\RHO}(\R^{2d})$.

 \section{Notation and Definitions.}
\label{sec: notation}
 In this section we introduce some notations and recall some standard definitions and state some main results from previous work.
\begin{itemize}
\item  $\SP(\R)$ is the set of all Borel probabilities on $\R.$

\item    $\SP_p(\R)$ ($1\leq p<\infty$) is the set of all probability measures $\mu\in\SP(\R)$ such that
$$\int_{\R} |x|^p d\mu<\infty.$$

  \item Given $\mu$ and $\nu \in \SP_p(\R)$, we denote by $\Gamma(\mu,\nu)$ the set of all Borel measures on $\R\times\R$  whose first and second marginal are respectively $\mu$ and $\nu$. 
  \item We say that a Borel map $\TF$ \textit{pushes} $\mu$ \textit{ forward} to  $\nu$   and write $\TF\#\mu=\nu$ if $\nu(A)=\mu(\TF^{-1}(A))$ for any Borel subset $A$ of $\R$. 
   \item Let  $\mu$, $\nu$ be two Borel measures of  finite $p-$moments. Then, the $p$-th Wasserstein distance between  $\mu$ and $\nu$  is defined by
   \begin{equation}\label{def of Wasserstein distance}
    \WAS^{p}_{p}(\mu,\nu)=     \min\left\lbrace \int_{\R^d\times\R^d}|x-y|^p d\gamma: \gamma\in \Gamma(\mu,\nu)\right\rbrace. 
   \end{equation}
  The set of minimizers in (\ref{def of Wasserstein distance}) is denoted by $\Gamma_0(\mu,\nu)$.
  \item A subset $\mathscr{S}\in \R\times\R$ is monotone if $\langle x_1-x_2, y_1-y_2\rangle \geq 0$ for any $(x_1,y_1), (x_2,y_2)\in \mathscr{S}. $
\end{itemize} 


\begin{definition}[Distance on $\SP_\RHO(\R)$]\label{def metric} 
Given $\RHO\in\SP_2(\R)$, we define 
$$
  \SP\SUB{\RHO}(\R^{2}) := \Big\{ \GAMMA\in\SP_2(\R^{2}) \colon
    \pi^1\#\GAMMA = \RHO \Big\}.
$$
The space $\SP_\RHO(\R)$ has a structure that looks like the one of Hilbert spaces  and is endowed with a distance in the following way : For any $\GAMMA^1,\GAMMA^2\in\SP
\SUB{\RHO}(\R^{2})$ with the disintegration $\GAMMA^k(dx,dy) =: \gamma^k_x(dy) \,\RHO(dx)$ with $k=1..2$, the distance $ \WAS\SUB{\RHO}(\GAMMA^1,\GAMMA^2)$ is defined by

$$
  \WAS\SUB{\RHO}(\GAMMA^1,\GAMMA^2)^2
    := \int_{\R^d} \WAS(\gamma^1_x,\gamma^2_x)^2 \,\RHO(dx).
$$

\end{definition}

\begin{definition}[Transport Plans]
Let $\RHO\in\SP_2(\R)$ be given.
\begin{enumerate}
\renewcommand{\labelenumi}{(\roman{enumi}.)}
\item {\em Admissible Plans.} For any $\GAMMA^1,\GAMMA^2\in\SP
  \SUB{\RHO}(\R^{2})$ we define
$$
  \ADM(\GAMMA^1,\GAMMA^2) := \Big\{ \ALPHA\in\SP(\R^{3}) \colon
    (\pi^1,\pi^2)\#\ALPHA = \GAMMA^1, 
    (\pi^1,\pi^3)\#\ALPHA = \GAMMA^2 \Big\}.
$$
\item {\em Optimal Plans.} For any $\GAMMA^1,\GAMMA^2\in\SP
    \SUB{\RHO}(\R^{2})$ we define
\begin{align*}
  \OPT(\GAMMA^1,\GAMMA^2) := \Big\{
    & \ALPHA\in\ADM(\GAMMA^1,\GAMMA^2) \colon
\\
    & \qquad \WAS\SUB{\RHO}(\GAMMA^1,\GAMMA^2)^2
      = \int_{\R^{3d}} |y^1-y^2|^2
        \,\ALPHA(dx,dy^1,dy^2) \Big\}.
\end{align*}
\end{enumerate}
\end{definition}

In the following theorem, we collect some fundamental results obtained in \cite{Gigli2004} about the distance in $\SP_\RHO(\R)$.

\begin{theorem}
Let $\RHO\in\SP_2(\R)$ be given.
\begin{enumerate}
\renewcommand{\labelenumi}{(\roman{enumi}.)}
\item The function $\WAS\SUB{\RHO}$ is a distance on $\SP\SUB{\RHO}(\R^{2})$
  and lower semicontinuous with respect to the weak convergence in
  $\SP(\R^{2})$. We have
$$
  \WAS\SUB{\RHO}(\GAMMA^1,\GAMMA^2)^2
    = \min_{\ALPHA\in\ADM(\GAMMA^1,\GAMMA^2)}
      \int_{\R^{3}} |y^1-y^2|^2 \,\ALPHA(dx,dy^1,dy^2)
$$
for all $\GAMMA^1, \GAMMA^2\in\SP \SUB{\RHO}(\R^{2})$, and thus
$\OPT(\GAMMA^1,\GAMMA^2)$ is nonempty.
\item The set $(\SP\SUB{\RHO}(\R^{2}), \WAS\SUB{\RHO})$ is a complete metric
  space.
\end{enumerate}
\end{theorem}


 The space $\SP_\RHO(\R)$ enjoys some vector space like properties when endowed with the following algebraic structure; see\cite{Gigli2004}.
\begin{definition}[Operations in $\SP\SUB{\RHO}(\R^{2})$]\label{def :algebraic structure} 
Let $\RHO\in\SP(\R)$ be given.
\begin{enumerate}
\renewcommand{\labelenumi}{(\roman{enumi}.)}
\item {\em Scalar multiplication.} For any $\GAMMA\in\SP\SUB{\RHO}(\R^{2})$ and $s\in\R$
  let
$$
  s\GAMMA := (\pi^1,s\pi^2)\#\GAMMA \in \SP\SUB{\RHO}(\R^{2}).
$$
\item {\em Addition.} For any pair $\GAMMA^1,\GAMMA^2\in\SP\SUB{\RHO} (\R)$ and $\BETA\in\ADM(\GAMMA^1,\GAMMA^2)$, we define the {\em addition} with respect to $\BETA$ by
$$  \GAMMA^1 +^{\BETA}\GAMMA^2 = (\pi^1,\pi^2+\pi^3)\#
    \BETA.    $$
  We denote the collection of all {\em additions} as $\BETA$ varies in $\ADM(\GAMMA^1,\GAMMA^2)$ by  
$$
  \GAMMA^1\oplus\GAMMA^2 := \Big\{ (\pi^1,\pi^2+\pi^3)\#
    \BETA\colon \BETA\in\ADM(\GAMMA^1,\GAMMA^2) \Big\}
      \subset\SP\SUB{\RHO}(\R).
$$
 Note that if either of $\GAMMA^i$ is induced by a map, then $\GAMMA^1\oplus\GAMMA^2 $ reduces to a singleton.
\end{enumerate}
\end{definition}

We introduce the following subset of $\SP\SUB{\RHO}(\R^{2})$ :
 
  \begin{equation}
  \COPT\SUB{\RHO} := \Big\{ \GAMMA \in \SP\SUB{\RHO}(\R^{2}) \colon 
    \text{$\SPT\GAMMA$ is a monotone subset of $\R\times\R$} \Big\}.
\label{E:CONE}
  \end{equation}
It is shown in \cite{CavSedWes2014} that $  \COPT\SUB{\RHO}$ is a closed convex cone with respect to the algebraic structure  and the distance above.
We recall the following theorem in \cite{Gigli2004} the proof of which uses the completeness of the metric space $(\SP\SUB{\RHO}(\R^{2}), \WAS\SUB{\RHO})$.  
  \begin{proposition}[Metric Projection]\label{P:metric proj}
Let $\RHO\in\SP_2(\R^d)$ be given and $C_{\RHO} \subset \PR$ a closed. For any $\gamma\in \PR$ there is a unique $\mathcal{P}_{\COPT\SUB{\RHO}} (\gamma)  \in \COPT\SUB{\RHO}$
with
\[
  \WR\big(\gamma,\mathcal{P}_{\COPT\SUB{\RHO}} (\gamma) \big) \LS \WR(\gamma,\ETA)
  \quad\text{for all $\ETA \in \COPT\SUB{\RHO}$.}
\]
\end{proposition}


\section{Preliminary.}\label{Prel}


 In this section, we consider the scaling map
 $\mathtt{S}^{\lambda}: x \longmapsto \lambda x$ and the translation map $\mathtt{T}^{u}: x\longmapsto x-u$ where $u\in \R$ and  $\lambda \in \R$. We state standard results on the actions of these maps on the Wasserstein distance.

\begin{lemma} \label{lem: scaling and tranlation}
Let $\mu, \nu\in \mathcal{P}_{2}(\R)$ and assume $\lambda> 0$. Then

$$W_{2}\left(\mathtt{S}^{\lambda}\#\mu,\mathtt{S}^{\lambda}\#\nu \right)=\lambda \WAS_{2}\left(\mu, \nu \right)\qquad\hbox{and}\qquad \WAS_{2}\left(\mathtt{T}^{u}\#\mu, \mathtt{T}^{u}\#\nu \right)=\WAS_{2}\left(\mu, \nu \right)$$
\end{lemma}
\proof
Let $\gamma\in\Gamma_0(\mu, \nu)$. Note that $(\mathtt{S}^{\lambda}, \mathtt{S}^{\lambda})\#\gamma\in \Gamma(\mathtt{S}^{\lambda}\#\mu, \mathtt{S}^{\lambda}\#\nu)$. Thus,

$$
\begin{aligned}
\int\vert x-y\vert^2 d(\mathtt{S}^{\lambda}, \mathtt{S}^{\lambda})\#\gamma= \int\vert\mathtt{S}^{\lambda}(x)-\mathtt{S}^{\lambda}(y)\vert^2d\gamma=  \int\vert \lambda x-\lambda y \vert^2d\gamma= \lambda^2 \WAS_{2}^2\left(\mu, \nu \right)
\end{aligned}$$
Thus, as  $(\mathtt{S}^{\lambda}, \mathtt{S}^{\lambda})\#\gamma\in \Gamma(\mathtt{S}^{\lambda}_\#\mu, \mathtt{S}^{\lambda}\#\nu),$
\begin{equation}\label{eq :scal1} 
\WAS_{2}^2\left(\mathtt{S}^{\lambda}\#\mu, \mathtt{S}^{\lambda}\#\nu \right)\leq \lambda^2 \WAS_{2}^2\left(\mu, \nu \right)
\end{equation}
By Symmetry, we have

\begin{equation}\label{eq :scal2}
\WAS^2\left(\mu, \nu \right)\leq \frac{1}{ \lambda^2 }\WAS^2\left(\mathtt{S}^{\lambda}\#\mu, \mathtt{S}^{\lambda}\#\nu \right) 
\end{equation}
We combine (\ref{eq :scal1}) and (\ref{eq :scal2}) to obtain  the first equality. We establish the second 
equality in a similar manner.   \endproof

\begin{lemma}\label{eq: W on trans. and scal.} 
 For $\lambda> 0$, $\gamma_{1}, \gamma_{2}\in \mathcal{P}_{\varrho}(\R)$, we have
  $$\WAS\SUB{\RHO} \left( \left( \pi^1,\lambda(\pi^2-\pi^1) \right)\# \GAMMA_1,\left(\pi^1,\lambda(\pi^2-\pi^1) \right)\#\GAMMA_2 \right) =\lambda \WAS\SUB{\RHO}(\GAMMA_1,\GAMMA_2)$$
  
  \end{lemma}
  
  \proof{}
 Note that  if $\GAMMA\in \SP_{\varrho}(\R)$ and $\bar\GAMMA=\left(\pi^1,\lambda(\pi^2-\pi^1) \right)\# \GAMMA$ then we have the following disintegration: 
  $$\bar\GAMMA=\int \bar\GAMMA_{x}d\varrho \quad\hbox{ with}\quad\bar\GAMMA_{x}=\lambda(\ID-x)\# \GAMMA_{x}$$
for $\varrho\;$ almost every $\, x$. In other words,
\begin{equation} \label{eq : trans comp scal}
\bar\GAMMA_{x}=S^{\lambda}\circ T^{x}\#\GAMMA_{x}.
\end{equation}
We use lemma \ref{lem: scaling and tranlation}  above to compute
 $$
\WAS^2\left(\mathtt{S}^{\lambda}\circ \mathtt{T}^{1x}\# \GAMMA_{1x},\mathtt{S}^{\lambda}\circ \mathtt{T}^{2x}\# \GAMMA_{2x} \right)=\lambda^2  \WAS^2\left(\GAMMA_{1x}, \GAMMA_{2x} \right).
$$ 
This, in view of (\ref{eq : trans comp scal}), yields
$$
  \begin{aligned}
  \WAS^2\SUB{\RHO}(\left(\pi^1,\lambda(\pi^2-\pi^1) \right)\# \GAMMA_1,\left(\pi^1,\lambda(\pi^2-\pi^1) \right)\#\GAMMA_2)
  = \lambda^2\WAS^2\SUB{\RHO}( \GAMMA_1,\GAMMA_2)
  \end{aligned}
$$ 
\endproof

\begin{proposition}
  Let $\lambda>0,\;\mu\in\SP\SUB{\RHO}(\R^{2})$. Set $\GAMMA_{\lambda}:=\left(\pi^{1},\pi^{1}+\lambda \pi^{2} \right)\#\mu $ and denote by $\mathcal{P}_{\COPT\SUB{\RHO}}(\GAMMA_{\lambda})$ its metric projection  onto $\COPT\SUB{\RHO}$ as provided by Proposition \ref{P:metric proj}. If ${\bf m}\in\COPT\SUB{\RHO}$ then
  
  \begin{equation}\label{E: inequality diff quotient} 
  \WAS\SUB{\RHO}\left( \left(\pi^1,\frac{\pi^2-\pi^1}{\lambda} \right)\#\mathcal{P}(\GAMMA_{\lambda}),\;\mu\right)\leq \WAS\SUB{\RHO}\left( \left(\pi^1,\frac{\pi^2-\pi^1}{\lambda} \right)\# {\bf m},\;\mu\right) 
  \end{equation}

  In particular, for $ \left\lbrace {\bf m}_{\lambda}\right\rbrace _{\lambda}\in\COPT\SUB{\RHO}$ 
  
  $$ \left(\pi^1,\frac{\pi^2-\pi^1}{\lambda} \right)\# {\bf m}_{\lambda}\xrightarrow { \WAS\SUB{\RHO}}\mu
   \implies\left(\pi^1,\frac{\pi^2-\pi^1}{\lambda} \right)\# \mathcal{P}(\GAMMA_{\lambda})\xrightarrow { \WAS\SUB{\RHO}}\mu $$
  
 \end{proposition}
 
 \proof
  By definition of the projection,
  \begin{equation}\label{E: inequality proj} 
  \WAS\SUB{\RHO}\left( \mathcal{P}(\GAMMA_{\lambda}),  \GAMMA_{\lambda}\right)\leq \WAS\SUB{\RHO}\left( {\bf m}, \GAMMA_{\lambda}\right)
  \end{equation}
 Applying lemma \ref{eq: W on trans. and scal.} to  (\ref{E: inequality proj}) we obtain (\ref{E: inequality diff quotient}).  
  \endproof
   Let  $\mu\in\SP_{\RHO}(\R^{2})$ and set 
  \begin{equation}
  \Lambda_{\mu} := \Big\{ \tau\in\R\colon \tau>0\; \text{and}\; \left(\pi^1,\pi^2+\tau\pi^1 \right)\#\mu \in  \COPT\SUB{\RHO} \Big\}.
  \end{equation}
 The next lemma states that either  $ \Lambda_{\mu} $ is empty or is an open interval.
 \begin{lemma}
 If $\tau_0\in  \Lambda_{\mu}$ and  $0<\tau <\tau_0$  then   $\tau \in\Lambda_{\mu}.$
 \end{lemma}
 
 \proof
Let $\tau, \tau_0> 0$ and set 

$$ \gamma_0= \left(\pi^1,\pi^2+\tau_0\pi^1 \right)\#\mu \quad \hbox{and} \quad\gamma= \left(\pi^1,\pi^2+\tau\pi^1 \right)\#\mu $$
 Note that, as  $\pi^1,\pi^2$ are continuous,
 $$ \text{supp}\gamma_0 =\overline{\left(\pi^1,\pi^2+\tau_0\pi^1 \right)\left[\text{supp}\mu\right]}. $$
  Let $(A_1,B_1)$, $(A_2,B_2)$  be  in the image of $ \text{supp} \mu$ by $\left(\pi^1,\pi^1+\tau\pi^2 \right)$ and choose $(x_1,y_1), (x_2,y_2)\in\text{supp} \mu$ such that
 
 $$ A_1=x_1,\quad A_2=x_2, \quad B_1=x_1+\tau y_1, \quad B_2=x_2+\tau y_2.$$
 We set,
 
 $$ R:=(A_2-A_1)(B_2-B_1)= (x_2-x_1)^2 + \tau (y_2-y_1)(x_2-x_1).$$ 
We claim that  $R\geq 0$ provided that $\tau_0\in  \Lambda_{\mu}$. To see this, assume  $$\tau_0\in  \Lambda_{\mu}\quad\hbox{and}\quad(y_2-y_1)(x_2-x_1) <0.$$ 
Thus,
  $$ R>(x_2-x_1)^2 + \tau_0 (y_2-y_1)(x_2-x_1) $$
 Therefore, as $(x_1, x_1+ \tau_0 y_1), (x_2, x_1+\tau y_2)$  belong to the image of $\text{supp}\mu$ by $ \left(\pi^1,\pi^1+\tau\pi^2 \right)$, the monotonicity of the support of $\left(\pi^1,\pi^2+\tau_0\pi^1 \right)\#\mu$ implies that
  
 \begin{equation}\label{eq: mon 4} 
 (A_2-A_1)(B_2-B_1)= R>(x_2-x_1)^2 + \tau_0 (y_2-y_1)(x_2-x_1) \geq 0
 \end{equation}
Since (\ref{eq: mon 4}) holds for any $(A_1,B_1)$, $(A_2,B_2)$  is in the image of $ \text{supp} \mu$ by $\left(\pi^1,\pi^1+\tau\pi^2 \right)$, we conclude that  the image of $ \text{supp} \mu$ by $\left(\pi^1,\pi^1+\tau\pi^2 \right)$, is monotone. As a result
$$  \quad \text{supp} \gamma =\overline{\left(\pi^1,\pi^2+\tau\pi^1 \right) \left[\text{supp}\mu \right]} $$
 is monotone by a density argument. It follows that $\gamma\in \COPT\SUB{\RHO}$ and so $\tau\in \Lambda_{\mu}.$

 
  \begin{lemma}\label{lem: density with control on supp} 
 Let $\mu\in\mathcal{P}_{\RHO}(\R^2)$. Then, there exists a sequence  $\lbrace \mu_n\rbrace_{n}\subset\SP_{\RHO}(\R^2)$ such that the support of $\mu_n$ in contained in $\R\times[-n,n]$ and, 
 
 $$ \mu_{n}\xrightarrow { \WAS\SUB{\RHO}}\mu
  $$
  \end{lemma}
  
\proof
 Write $\mu=\int \mu^{x}d\rho$ and  define $ P_{n}(\xi):=\xi\chi_{\vert \xi \vert \leq n}$.
Consider a sequence of measures $\lbrace \mu_n\rbrace_{n}\subset\mathcal{P}_{\RHO}(\R^2)$ such that  $\mu_n^{x}=P_n\#\mu^{x}$. Clearly, the support of $\mu_n$ in contained in $\R\times[-n,n]$ and
 
$$W^2( \mu_n^{x}, \mu^{x} )\leq \int_{\R}\vert P_n(\xi)-\xi \vert^{2} d\mu^{x}=\int_{\vert\xi\vert>n}\vert \xi \vert^{2} d\mu^{x} $$

As, $x\longrightarrow \int_{\vert\xi\vert>n}\vert \xi \vert^{2} d\mu^{x} \in L^1(\varrho)$ and converges to $0$, by the Lebesgue dominated convergence theorem, $W_{\RHO}( \mu_n, \mu )$ converges to $0$.\endproof


 In the sequel, let us write $\RHO\in\SP_{2}(\R^2)$ as
 $$\RHO=\sum^{\infty}_{i=1} m_i\delta_{x_i}+\bar{\RHO}.$$
 where $\left\lbrace x_i\right\rbrace^{\infty}_{i=1}\subset \R $ is the set of atoms of $\RHO$ and $\bar{\RHO}$ is the diffuse part of $\RHO$ and $ m_i \geq 0$ so that 
 $$ \sum_{i=1}^{\infty} m_i+ \bar\RHO (\R)=1. $$
 
\begin{lemma}\label{lem: density with control of velocity at atoms} 
 Let $\mu\in\SP_{\RHO}(\R^2)$ such that 
\begin{equation}\label{E: standard mu} 
\mu=\sum^{\infty}_{i=1}\nu_{i}\times\delta_{x_i}+ (\ID, u)\#\bar{\RHO}.
\end{equation}
for some $u\in L^2(\bar{\RHO})$ and Borel measures $\nu_i$ on $\R$.
Set 
$$\mu_{N}=\sum^{N}_{i=1}\nu_{i}(\xi)\times\delta_{x_i}+\sum^{\infty}_{i=N+1}\delta_{0}(\xi)\times\delta_{x_i}+ (\ID, u)\#\bar{\RHO}. $$
 Then,
 
 $$ \mu_{N}\xrightarrow { \WAS\SUB{\RHO}}\mu
  $$
 
  \end{lemma}
\proof  Note that 
\begin{equation}\label{E: conv series W_2 1} 
\WAS^2( \delta_{0}, \nu_{i} )=\int_{\R}\vert \xi \vert^{2} d\nu_{i}
\end{equation}
 Using (\ref{E: standard mu}), we obtain that
\begin{equation}\label{E: conv series W_2 2} 
\sum^{\infty}_{i=1} \int_{\R}\vert \xi \vert^{2} d\nu_{i}\leq  \int_{\R}\vert \xi \vert^{2} d\mu
\end{equation}
 Since
 \begin{equation}\label{E: conv series W_2 3} 
 \WAS^2\SUB{\RHO}( \mu_N, \mu )=\sum^{\infty}_{i=N+1}  \WAS^2( \delta_{0}, \mu^{x_i} )= \sum^{\infty}_{i=N+1}  \WAS^2( \delta_{0}, \nu_{i}).
\end{equation}
As $\mu$ is of finite second moment, we combine (\ref{E: conv series W_2 1})-(\ref{E: conv series W_2 3}), to obtain that 
$$\lim_{N\longrightarrow\infty}  \WAS^2\SUB{\RHO}( \mu_N, \mu )= 0. $$
\endproof

In order to state our main result, we introduce the tangent space of the monotone transport plans: 
 
 $$ \mathbb{T}_{\COPT\SUB{\RHO}}:= \overline{\left\lbrace \mu\in \SP_{\varrho}(\R^2) : (\pi^1,\pi^1+\tau \pi^2)\#\mu \in\COPT\SUB{\RHO}\; \text{ for some } \tau>0 \right\rbrace}^{  \WAS\SUB{\RHO}}  $$

 \begin{theorem}
 $ \mathbb{T}_{\COPT\SUB{\RHO}}$ is a closed convex  cone with respect to the algebraic structure given in Definition \ref{def :algebraic structure} . Furthermore, $\mu\in  \mathbb{T}_{\COPT\SUB{\RHO}}$ if and only if $\mu$ has the following representation :
  
 \begin{equation}\label{eq: linking map} 
 \mu = (\ID,g)\#\bar{\RHO} + \sum^{\infty}_{i=1}( \nu_i\times\delta_{y_i})
 \end{equation}
  for some $g\in L^2(\bar{\RHO})$, a sequence $\left\lbrace y_i\right\rbrace_{i=1}^{\infty}\subset \R $ and some family of Borel probability measures $\left\lbrace \nu_i\right\rbrace_{i=1}^{\infty} $ on $\R$.
  
 \end{theorem}
 The proof of this theorem follows:


\section{ The tangent cone of monotone transport plans is  convex.}\label{Useful}

In this section we prove the part of the main theorem that is concerned with the fact the tangent space is a convex cone.
   \begin{lemma}\label{lem: gigli}
 Let $\mu_i, \bar\mu_i\in\SP\SUB{\RHO}(\R^{2})$ $i=1,2$ and $\mu\in \mu^1\oplus\mu^2 $. Then there exists $\bar\mu\in \bar\mu^1\oplus\bar\mu^2$ such that 
 
 $$\WAS\SUB{\RHO}(\mu,\bar\mu) \leq\WAS\SUB{\RHO}(\bar\mu^2,\mu^2) +\WAS\SUB{\RHO}(\mu^1,\bar\mu^1) $$
 \end{lemma}
 \proof we refer the reader to (\cite{Gigli2004} Proposition 4. 24).
\endproof

 \begin{proposition}\label{Prop: convex cone} 
 \begin{enumerate}
 \item $ \mathbb{T}_{\COPT\SUB{\RHO}}$ is a  cone.
 \item $ \mathbb{T}_{\COPT\SUB{\RHO}}$ is convex : if $\mu_1, \mu_2\in \mathbb{T}_{\COPT\SUB{\RHO}}$ then $\mu^1\oplus\mu^2 \subset  \mathbb{T}_{\COPT\SUB{\RHO}}$.
 \end{enumerate}
 .
 \end{proposition}
 
  \proof
(i)  $\lambda\in\mathbb{\R}_{+}$ and $\mu\in \mathbb{T}_{\COPT\SUB{\RHO}}$. Then, there exist $\left\lbrace \mu_{n} \right\rbrace_{n=1}^{\infty}  $ and $\left\lbrace \tau_{n} \right\rbrace_{n=1}^{\infty}$ such that $\left\lbrace \mu_{n} \right\rbrace_{n=1}^{\infty}  $  converges to $\mu$ with respect to $\WAS\SUB{\RHO}$ and $\tau_{n}\in \Lambda_{\mu_{n}}$. Set $\nu_{n}=(\pi^1,\lambda\pi^2)\#\mu_{n}$ and note that 

$$\left( \pi^1,\pi^1 +\frac{\tau_n}{\lambda}\pi^2\right)\#\nu_{n} = \left( \pi^1,\pi^1 +\tau_n\pi^2\right)\#\mu_{n}\in \COPT\SUB{\RHO}$$
  Thus,  $\frac{\tau_n}{\lambda}\in\Lambda_{\nu_{n}}$. On the other hand, using lemma \ref{lem: scaling and tranlation}, we easily show that
 
$$\WAS\SUB{\RHO}(\nu_n,\nu) =\lambda \WAS\SUB{\RHO}(\mu_n,\mu) $$ 
  with $\nu:= (\pi^1,\lambda\pi^2)\#\mu$. Thus, $\left\lbrace \nu_{n} \right\rbrace_{n=1}^{\infty}$ converges to $\nu$ with respect to $\WAS\SUB{\RHO}$. So, $\lambda \mu\in\mathbb{T}_{\COPT\SUB{\RHO}}$. we conclude that  $ \mathbb{T}_{\COPT\SUB{\RHO}}$ is a cone.\\
  (ii)
$\mu^i\in\mathbb{T}_{\COPT\SUB{\RHO}}$  $i=1,2$. There exists  a positive sequence $\left\lbrace  \tau_{n}\right\rbrace^{\infty}_{n=1}$ and   $\left\lbrace  \mu^i_{n}\right\rbrace^{\infty}_{n=1}\subset\in\SP
  \SUB{\RHO}(\R^{2})
  $ $i=1,2$ such that $\left\lbrace  \mu^i_{n}\right\rbrace^{\infty}_{n=1}$ converges to  $\mu^i$ with respect to $\WAS\SUB{\RHO}$ and 

$$ \GAMMA^i_{n}=\left( \pi^1, \pi^1+\tau_{n}\pi^2\right)\#\mu^i_{n}\in \COPT\SUB{\RHO}.$$
Let $\mu\in \mu^1\oplus\mu^2$. Then, by lemma \ref{lem: gigli}, there exists $\mu_n\in \mu_n^1\oplus\mu_n^2$
such that 
 \begin{equation}\label{eq: conv}
  \mu_{n}\xrightarrow { \WAS\SUB{\RHO}}\mu.
  \end{equation}
Let $\beta_{n}\in \ADM(\mu_n^1,\mu_n^2)$ so that 

$$\mu^1_{n}=(\pi^1,\pi^2)\#\beta_{n},\; \mu^2_{n}=(\pi^1,\pi^3)\#\beta_{n},\; \hbox{ and } \; \mu^1_{n}+^{\beta_n}\mu^2_{n}=(\pi^1,\pi^2+\pi^3)\#\beta_{n}, $$
We use  the identities $$\left( \pi^1, \pi^1+\tau_{n}\pi^2\right)\circ(\pi^1,\pi^2)=\left( \pi^1, \pi^1+\tau_{n}\pi^2\right)$$
 and $$ \left( \pi^1, \pi^1+\tau_{n}\pi^3\right)\circ(\pi^1,\pi^2)=\left( \pi^1, \pi^1+\tau_{n}\pi^3\right)$$
  to get respectively 
  $$\GAMMA^1_{n}=\left( \pi^1, \pi^1+\tau_{n}\pi^2\right)\#\beta_{n} \quad \hbox{ and } \quad \GAMMA^2_{n}=\left( \pi^1, \pi^1+\tau_{n}\pi^3\right)\#\beta_{n}   $$
Note that 
$$ \begin{aligned}
\left( \pi^1, \pi^1+\frac{\tau_n}{2}\pi^2\right)\#\mu_n &=\left( \pi^1, \pi^1+\frac{\tau_n}{2}\pi^2\right)\#(\pi^1,\pi^2+\pi^3)\#\beta_{n}\\
&=\left( \pi^1, \pi^1+\frac{\tau_n}{2}(\pi^2+\pi^3)\right)\#\beta_n\\
\end{aligned}
$$
That is,
\begin{equation}\label{eq: 2nd cond for Tan}
 \begin{aligned}
 \left( \pi^1, \pi^1+\frac{\tau_n}{2}\pi^2\right)\#\mu_n & = \frac{1}{2} \left( \pi^1, \pi^1+\tau_{n}\pi^2\right)\#\beta_n + \frac{1}{2} \left( \pi^1, \pi^1+\tau_{n}\pi^3\right)\#\beta_n \\
&= \frac{1}{2}(\GAMMA^1_{n}+\GAMMA^2_{n})
\end{aligned}
   \end{equation}
As, $ \GAMMA^i_{n}\in \COPT\SUB{\RHO}$ and $\COPT\SUB{\RHO}$ is convex,  (\ref{eq: 2nd cond for Tan}) implies that 
$$ \left( \pi^1, \pi^1+\frac{\tau_n}{2}\pi^2\right)\#\mu_n \in \COPT\SUB{\RHO}. $$
This combined with (\ref{eq: conv}) yields that $\mu\in \mathbb{T}_{\COPT\SUB{\RHO}}$. And so, $\mu^1\oplus\mu^2\subset \mathbb{T}_{\COPT\SUB{\RHO}}.$ \endproof

  
    
    
     
     

\section{Characterization of  the Tangent space.}\label{Chara of the Tangent} 
In this section, we give a complete characterization of the tangent space of monotone plans.  Using the convexity property  established in the previous section, we show that elements of the tangent cone of monotone plans are made essentially of two basic components. 


  We recall
  
 $$\RHO=\sum^{\infty}_{i=1} m_i\delta_{x_i}+\bar{\RHO}.$$
 where $\left\lbrace x_i\right\rbrace^{\infty}_{i=1}\subset \R $ is the set of atoms of $\RHO$ and $\bar{\RHO}$ is the diffuse part of $\RHO$ and $m_i \geq 0$ so that 
 $$ \sum_{i=1}^{\infty} m_i+ \bar\RHO (\R)=1. $$
 
  \begin{lemma}\label{lem : measures in Tan} 
 Let $\gamma\in\COPT\SUB{\RHO}$. Then,
  there exists monotone function $u$, a countable set $B :=\left\lbrace y_i \right\rbrace^{\infty}_{i=1} \subset \R$  and a  family of Borel probability measures  $\left\lbrace \nu_i \right\rbrace^{\infty}_{i=1} $  such that 
  \begin{equation}\label{eq: representation of gamma monotone}
  u\in L^2 (\bar{\RHO}) \quad \hbox{and}\quad \GAMMA= (\ID, u)\#\bar \RHO+\sum^{\infty}_{i=1}\nu_i\times\delta_{y_i}. 
  \end{equation}

  \end{lemma}
\proof
 Let $\Gamma$ be a maximal monotone extending set  for $spt (\GAMMA)$ and $u$ the corresponding set valued map.  Let $A$ be the set of all points $x$ at which $u(x)$ has more than one point.  As $u(x)$ is  a closed interval in $\R$,  $A$ has at most countably many points, that is, $A=\left\lbrace m_i \right\rbrace^{\infty}_{i=1}$. Let $K$ be a $\RHO$ measurable set such that $\RHO(K^c)=0$ and the disintegration of $\GAMMA$ with respect to $\RHO$ is uniquely given as $\GAMMA = \int \GAMMA_{x} d\RHO$ with $\left\lbrace \GAMMA_{x} \right\rbrace_{x\in K}.$  Let $x\in A^{c}$ so that $u(x)$ is a singleton. If, in addition, $x\in K$ then $\GAMMA_{x}=\delta_{u(x)}$.
  Clearly $$\GAMMA =\int \delta_{u(x)}d \RHO + \sum_{i=1}^{\infty} \GAMMA_i\times \delta_{m_i}$$
  for some family of Borel probability measures $\left\lbrace \gamma_i \right\rbrace^{\infty}_{i=1}.$
  As $\RHO=\sum^{\infty}_{i=1}m_i\delta_{x_i}+\bar{\RHO}$, we have
  
  $$\GAMMA= \int \delta_{u(x)}d \bar\RHO +\sum^{\infty}_{i=1}\nu_i\times\delta_{y_i} $$
  for some $\left\lbrace y_i \right\rbrace^{\infty}_{i=1} \subset \R$. The fact that $\GAMMA$ has a finite second moment ensures that $u \in L^2(\bar\RHO)$.
  
  \endproof

\begin{proposition} \label{lem: representation of mu}
 If $\mu\in \mathbb{T}_{\COPT\SUB{\RHO}}$ then $\mu$  has the following representation :

 \begin{equation}\label{eq: linking map} 
 \mu = (\ID,g)\#\bar{\RHO} + \sum^{\infty}_{i=1}( \nu_i\times\delta_{y_i})
 \end{equation}
  for some $g\in L^2(\bar{\RHO})$ and family of Borel probability measures $\left\lbrace \nu_i\right\rbrace^{\infty}_{i=1} $ on $\R$.
  
 \end{proposition}
 
 \begin{remark}
 We are going to show below that any  $g\in L^2(\bar{\RHO})$ and any choice of family of Borel probability measures $\left\lbrace \nu_i\right\rbrace^{\infty}_{i=1} $ will give rise to elements of  $ \mathbb{T}_{\COPT\SUB{\RHO}}$
 \end{remark}

\proof
 Let  $\mu\in \mathbb{T}_{\COPT\SUB{\RHO}}$. Then, there exists $\left\lbrace \mu_{n}\right\rbrace _{n=1}^{\infty}\subset  \SP\SUB{\RHO}(\R^{2})$ converging to $\mu$ and  a sequence of positive numbers $\left\lbrace \tau_{n} \right\rbrace^{\infty}_{i=n}$
  such that 
  $$\GAMMA_n:=(\pi^1,\pi^1+\tau_n \pi^2)\#\mu_n \in \COPT\SUB{\RHO}. $$
By lemma \ref{lem : measures in Tan},$\GAMMA_n$ have the  representation (\ref{eq: representation of gamma monotone}). Since $\mu_n=\left(\pi^1, \frac{\pi^2-\pi^1}{\tau_n} \right)\#\gamma_n $, $\mu_n$  have the same representation as (\ref{eq: representation of gamma monotone}) without $u$ being necessary monotone. Thus, there exists a family of functions $\left\lbrace u_n \right\rbrace^{\infty}_{i=1}\subset  L^2(\bar{\RHO})$,   $ A_n=\left\lbrace y_i^n \right\rbrace^{\infty}_{i=1}\subset \R$ and  a family of probability measures $\left\lbrace \nu^n_i \right\rbrace^{\infty}_{i=1}$ such that 
  $$\mu_n= \int_{\R} \delta_{u_n(x)}d \bar\RHO +\sum^{\infty}_{i=1}\nu^n_i\times\delta_{y^n_i}. $$
By setting $A=\cup A_n$, we rewrite $\GAMMA_n$ as
 
  $$\mu_n= \int_{A^c} \delta_{u_n(x)}d \bar\RHO +\sum^{\infty}_{i=1}\bar\nu^n_i\times\delta_{s_i} $$
    where $ A=\left\lbrace s_i \right\rbrace^{\infty}_{i=1}$ and  $\bar\nu^n_j=\nu^n_i$ if  $s_j=y^n_i$.
    
  Note that
  $$  \WAS^2\SUB{\RHO}( \mu_n, \mu_m )=\int_{A^c}\vert u_n-u_m \vert^{2} d\bar\RHO +  \sum^{\infty}_{j=1} W^2(\bar\nu^n_j,\bar\nu^m_j). $$
As $\left\lbrace \mu_n \right\rbrace^{\infty}_{i=1}$ is Cauchy, using completeness of the $L^2(\bar\RHO)$ and $\mathcal{P}_{2}(\R)$, we obtain that $\left\lbrace \mu_n \right\rbrace^{\infty}_{i=1}$  converges to

\begin{equation} 
 \mu = (\ID,g)\#\bar{\RHO} + \sum^{\infty}_{i=1}( \nu_i\times\delta_{s_i})
 \end{equation}
  for some $g\in L^2(\bar{\RHO})$ and borel measures $\nu_i$ on $\R$. The uniqueness of the limit ensures the result.
\endproof
   
 \begin{lemma}\label{lem : function in Tan} 
 Let $ g\in L^2(\RHO)$. Then,
 
 \begin{equation}\label{eq: function in Tanspace} 
 \mu:=(\ID,g)\#\RHO \in\mathbb{T}_{\COPT\SUB{\RHO}}
 \end{equation}

  \end{lemma}
\proof

There exists $\left\lbrace g_n\right\rbrace^{\infty}_{n=1}\subset \mathcal{C}_{c}^{\infty}(\R) $  such that $\left\lbrace g_n\right\rbrace^{\infty}_{n=1}$ converges to $g$ in $L^2(\RHO)$. Denoting $\mu_{g_n}:=(\ID,g_n)\#\RHO$,  we note that 
$$\WAS\SUB{\RHO}^2\left(\mu_{g_n},\mu\right)= \int_{\R} \vert g-g_n\vert^2 d\RHO.$$
Thus, $\left\lbrace \mu_{g_n}\right\rbrace^{\infty}_{n=1}$ converges to $\mu$ with respect to $\WAS\SUB{\RHO}$. As $\left\lbrace g_n\right\rbrace^{\infty}_{n=1}\subset \mathcal{C}_{c}^{\infty}(\R) $, we choose a sequence of  positive numbers $\left\lbrace \tau_n\right\rbrace^{\infty}_{n=1}$ such that  $1+\tau_n \Vert g'_n\Vert_{\infty} >0$. Then, we obtain that 
 
$$\gamma_n=(\pi^1,\pi^1+\tau_n \pi^2)\#\mu_{g_n}\in\COPT\SUB{\RHO}. $$
 We conclude that $\mu\in \mathbb{T}_{\COPT\SUB{\RHO}}.$ \endproof
 

    \begin{lemma}\label{lem : simple measure in Tan} 
Let $x_0$ be an atom of $\RHO$ and write $\RHO= m_0\delta_{x_0}+ \tilde{\RHO}$ so that $m_0+ \tilde{\RHO}(\R)=1$. Set 
 \begin{equation}\label{eq:  simple measure in Tan} 
   \mu:= m_0\nu\times\delta_{x_0} +\delta_{0}\times \tilde{\RHO}
  \end{equation}
where $\mu$ is a Borel probability measure with compact support. Then $  \mu\in\mathbb{T}_{\COPT\SUB{\RHO}}$.

  \end{lemma}
  
\proof

Let $\alpha>0$ such that $[-\alpha,  \alpha]$ is the smallest symmetric  interval  containing the support of $\nu$. We consider an increasing sequence $\left\lbrace \alpha_n \right\rbrace_{n=1}^{\infty} \subset (0,\alpha)$ such that $\alpha_n $ converges to $\alpha$. Define 

$$
S_{n}(x):=
\begin{cases}
x\quad \hbox{on}\quad (-\infty,  x_0-\frac{\alpha_n}{n}]\cup  [ x_0+\frac{\alpha_n}{n}\infty)\\
x_0-\frac{\alpha_n}{n} \qquad \hbox{on}\qquad  (x_0-\frac{\alpha_n}{n}, x_{0})\\
x_0+\frac{\alpha_n}{n}\quad \hbox{on}\quad  [x_0, x_0+\frac{\alpha_n}{n})\\
\end{cases}
$$
 and
 $$
P_{n}(\xi):=
\begin{cases}
x_0-\frac{\alpha_n}{n} \quad \hbox{on}\quad (-\infty, -\alpha_n)\\
x_0+\frac{\xi}{n}\qquad \hbox{on}\qquad  [-\alpha_n,\alpha_n ]\\
x_0+\frac{\alpha_n}{n} \quad \hbox{on}\quad  (\alpha_n ,\infty).\\
\end{cases}
$$
  We note that $S_n$ is monotone increasing with exactly one jump at $x_0$ of gap the interval $[x_0-\frac{\alpha_n}{n}, x_0+\frac{\alpha_n}{n}],$. We use $P_n$, $S_n$  to construct the following measure :
 $${\bf m}_{n}:= m_0 P_{n\#}\nu(d\xi)\times \delta_{x_0}(dx)+ \left(\ID,S_{n}\right)\# \tilde\RHO(dx).$$
 As $S_n$  is monotone and $P_n$ has values in $[x_0-\frac{\alpha_n}{n}, x_0+\frac{\alpha_n}{n}],$  it is immediate that the support of $ {\bf m}_{n}$ is monotone.\\
We derive a truncation for $\mu$, by setting 
 $$\mu_{n}(dx, d\xi)=\left(\pi^{1}, + \frac{\pi^{2}-\pi^{1}}{\frac{1}{n}}\right)\#{\bf m}_{n}.$$ 
 That is,
 $$\mu_{n}(dx, d\xi)=\sum^{\infty}_{i=1} m_i\delta_{\bar S_n(x_i)}(d\xi)\times\delta_{x_i}(dx) +\bar P_{n\#}\nu(d\xi)\times \delta_{x_i}(dx)+ \left(\ID,\bar S_n\right)\# \bar\rho(dx).$$
Here, 
$$\bar S_n= \frac{S_n-\ID}{\frac{1}{n}}\qquad \hbox{and} \qquad \bar P_n= \frac{P_n-\ID}{\frac{1}{n}}. $$
 A simple computation gives
$$
\bar S_{n}(x):=
\begin{cases}
0\quad \hbox{on}\quad (-\infty,  x_0-\frac{\alpha_n}{n})\cup  ( x_0+\frac{\alpha_n}{n}, \infty)\\
\frac{x_0-x}{\frac{1}{n}}-\alpha_n \qquad \hbox{on}\qquad  [x_0-\frac{\alpha_n}{n}, x_{0}]\\
\frac{x_0-x}{\frac{1}{n}}+\alpha_n \quad \hbox{on}\quad  (x_0, x_0+\frac{\alpha_n}{n}].\\
\end{cases}
$$
and 
$$
\bar P_{n}(\xi):=
\begin{cases}
-\alpha_n\quad \hbox{on}\quad (-\infty, \alpha_n)\\
\xi\qquad \hbox{on}\qquad  \left[-\alpha_n, \alpha_n\right]\\
\alpha_n\quad \hbox{on}\quad (\alpha_n, \infty).\\
\end{cases}
$$
Clearly, as $\alpha_n$ converges to $\alpha$,
 \begin{equation}\label{eq: limit nu part} 
\limsup_{n\rightarrow \infty} W^2\left(\bar P_{n\#}\nu,\nu\right))\leq \limsup_{n\rightarrow\infty}\int_{\R}\vert\bar P_{n\#}(\xi)-\xi\vert^{2}d\nu=0.
 \end{equation}
Observe that 

 \begin{equation}\label{eq: bound on Bar S} 
\vert\bar S_{n}(x)\vert\leq \vert\frac{x_0-x}{\frac{1}{n}}\vert+\alpha_n \leq 2\alpha
\end{equation}
and

 \begin{equation}\label{eq: lim Sn} 
\lim_{n\rightarrow\infty} \bar S_{n}(x)=0 \qquad \tilde{\RHO}- a.e.
 \end{equation}
 By using the Lebesgue dominated convergence theorem, (\ref{eq: bound on Bar S}) and (\ref{eq: lim Sn}) lead to
  \begin{equation}\label{eq: limit integral part} 
 \lim_{n\rightarrow\infty}\int_{\R} \mid\bar S_{n}(x)\mid^2 d\bar\RHO=0.
 \end{equation}
Since
  \begin{equation}\label{eq: W_rho into 2 parts}
   \WAS\SUB{\RHO}^2\left(\mu_{n},\mu\right)= \int_{\R} \mid\bar S_{n}(x)\mid^2 d\tilde\RHO + m_0 W^2\left(\bar P_{n\#}\nu,\nu\right),  
  \end{equation}
We use  (\ref{eq: limit nu part}) and (\ref{eq: limit integral part}) to obtain

 $$ \mu_{n}\xrightarrow { \WAS\SUB{\RHO}}\mu.
  $$ \endproof
 Now, let us  give a proof of the main theorem stated in section \ref{Prel}.\\
 \textit{Proof of the main theorem.}\\
The fact that $\mathbb{T}_{\COPT\SUB{\RHO}}$ is a convex cone is established in Proposition \ref{Prop: convex cone}. In view of the results in Proposition \ref{lem: representation of mu}, we only need to show the converse statement. Assume $\mu$ has the representation (\ref{eq: linking map}). By lemma \ref{lem: density with control on supp}  and lemma \ref{lem: density with control of velocity at atoms}, it suffices to show that 

 \begin{equation}\label{eq: finite sum of measure in Tan} 
 \mu := (\ID,g)\#\bar{\RHO} + \sum^{N}_{i=1}( \nu_i\times\delta_{x_i})+\sum^{\infty}_{i=N+1}( \delta_{0}\times\delta_{x_i})
 \end{equation}
  for some   positive integer $N$, $g\in L^2(\bar{\RHO})$ and Borel measures $\nu_i$ with support contained in  closed and bounded intervals, is an element of $\mathbb{T}_{\COPT\SUB{\RHO}}$. Furthermore, as  $\mathbb{T}_{\COPT\SUB{\RHO}}$ is a convex cone, the proof reduces to showing that $\mu$ in (\ref{eq: function in Tanspace}) and (\ref{eq:  simple measure in Tan}) belong to  $\mathbb{T}_{\COPT\SUB{\RHO}}$. These are obtained  respectively in lemma \ref{lem : function in Tan}  and lemma \ref{lem : simple measure in Tan}.  \endproof


\end{document}